\documentclass[12pt]{article}
\usepackage{mathrsfs}
\usepackage{amsmath}
\usepackage{amsfonts}
\usepackage{amssymb, amsmath, cite}
\usepackage{bm}
 \usepackage{verbatim}
\usepackage{amsthm}

\newtheorem{theorem}{Theorem}
\newtheorem{lemma}{\hskip\parindent\bf{Lemma}}[section]
\newtheorem{proposition}{Proposition}
\newtheorem{remark}{Remark}

\renewcommand{\theequation}{\thesection.\arabic{equation}}

\newcommand\be{\begin{equation}}
\newcommand\ee{\end{equation}}
\newcommand\ber{\begin{eqnarray}}
\newcommand\eer{\end{eqnarray}}
\newcommand\berr{\begin{eqnarray*}}
\newcommand\eerr{\end{eqnarray*}}
{\setlength\arraycolsep{2pt}

\begin{document}

\title{ \bf\Large  Existence of global vortices in a class of Ginzburg-Landau models }
\date{}
{\tiny\author{
  Xiaosen Han$^1$\ \ \ \ \   Kexin Zhang$^2$
\\
\\
 $^1$
Institute of Contemporary Mathematics, \\School of Mathematics and Statistics,\\ Henan University
Kaifeng 475004, PR China \\
$^2$Chern Institute of Mathematics,\\ Nankai University Tianjin 300071, PR China }
\maketitle}

\vspace{0.1in}
\bigskip
\begin{quote}
{{{\bfseries Abstract.}  In this paper, the two-component Ginzburg-Landau (TCGL) model, which is an important research object in mathematics and physics of two scalar fields with $U(1)\times U(1)$ symmetry, is studied in detail by the shooting method and the fixed point theorem, which contains two cases: single vacuum expectation value (1VEV) and two vacuum expectation values (2VEV). This represents two different physical states: superconductivity-normal state and superconductivity-superconductivity state, and we prove the existence, uniqueness, asymptotic properties and quantization of the solutions to these two boundary value problems.
}

}

\end{quote}

\section{Introduction}\label{s1}

The Ginzburg-Landau model \cite{ginzburg1950} is an important model describing superconductivity theory after Onnes' discovered superconductors \cite{onnes1911}.
In the mid-20th century, through the energy gap problem, Bardeen, Cooper and Schrieffer (BCS) \cite{bardeen1957} proposed that the pairing of two electrons, with opposite spin matching between arbitrary small attraction is the source of superconductivity by using quantum mechanics. The condensate is made of Cooper pairs, which can be equivalently described as a composite boson, similar to Bose-Einstein condensation\cite{kasamatsu2005,kasamatsu2016}. This is called BCS theory for type \uppercase\expandafter{\romannumeral1} superconductor, which has been confirmed by many experiments \cite{anderson1958,anderson1959,nambu1961,arkady2001}. For type \uppercase\expandafter{\romannumeral1} superconductor, if the magnetic field strength is low enough, the magnetic field will be excluded from the material, which is called the Meissner effect.  Almost at the same time, Abrikosov \cite{abrikosov1957} brought forward that the magnetic properties of superconducting alloys and predicted that when the external magnetic field strength exceeds the critical magnetic field strength, there is periodic structure with the symmetry of a 2-dimensional lattice and the interface energy would be negative between the superconductor and a normal metal. Also, there are many developments \cite{takayama1973} and applications \cite{chou1958,phillips1959} for type \uppercase\expandafter{\romannumeral2} superconductor. It may form a 'mixed state' in which the magnetic field penetrates the metal partially.  It is quite different from the properties of the type \uppercase\expandafter{\romannumeral1} superconductor, and their differences are due to the special vortex structure formed by type \uppercase\expandafter{\romannumeral2} superconductor, which plays an important part in physics and mathematics \cite{berlinsky1995,lin2000,lin2002,zhitomirsky2004}. The core of the vortex is defined by the superconductor's coherence length $\xi$ and London penetration depth $\lambda$, both of which are functions of temperature. Apart from the difference in microstructure, the Ginzburg-Landau parameter $\kappa=\lambda/\xi$ can be used to distinguish the two superconductors. In general, $\kappa$ is small for pure metals, but large for superconducting alloys. When $\kappa<1/\sqrt{2}$, it is of type \uppercase\expandafter{\romannumeral1}, and the interaction between the vortices is attractive. When $\kappa=1/\sqrt{2}$, it is called Bogomol 'nyi point \cite{Annett2004,Forgacs2006,forgacs2006} or critical point between type \uppercase\expandafter{\romannumeral1} and \uppercase\expandafter{\romannumeral2}. In this case, there is no interaction between the vortices. When $\kappa>1/\sqrt{2}$, it belongs to type \uppercase\expandafter{\romannumeral2}, and the interaction between the vortices is repulsive.

At the same time, there are many mathematical results about the classical one-parameter Ginzburg-Landau equation
\begin{equation}
-\Delta u+(|u|^2-1)u=0.\label{0}
\end{equation}
For example, the unique radial solution for \eqref{0} with Dirichlet boundary conditions in a simply connected region of $\mathbb{R}^2$, by variation method and shooting method, has been obtained in \cite{bookbrezis} and \cite{chen1994}, respectively. Quantization effects in $\mathbb{R}^2$ for \eqref{0} are given in \cite{brezis1994}. Liouville-type theorem and symmetric on entire solution and $O(N)$-equivariant vortex solution of \eqref{0} in higher dimensions has been obtained in \cite{farina2004,pisante2011}. Taubes answered Weinberg's conjecture by proving existence and uniquness of arbitrary $N$-vortex solutions to the one parameter Ginzburg-Landau equations \cite{taubes1980}.

With the discovery of liquid metal hydrogen, magnesium diboride, iron base, etc., superconductors (multi-band superconductors) with more than one condensate have come into our sight. Subsequently, the multi-component Ginzburg-Landau (MCGL) theory was proposed, providing a description for multi-band superconductors. MCGL theory has been applied to high-energy physics, condensed matter systems and astrophysics, providing a method for observing and studying astrophysical fields, plasma, and even the modeling of neutron stars \cite{jones2006}. We focus on the TCGL system which is already interesting in many applications.
Also vortices in TCGL system could carry fractional magnetic flux \cite{babaev2002}.
Unlike the other two types of superconductors, these kind of materials require higher phase transition temperatures and pressures, and are known as type 1.5 superconductivity \cite{babaev2004,moshchalkov2009,babaev2012}. A famous example of this type of superconductor is liquid metal hydrogen, which converts the phase into a superconducting state under high pressure. However, its electrons and protons form two different types of Cooper pairs at this point. To characterize this phenomenon, current-carrying vortices of the TCGL model with a $U(1)\times U(1)$ symmetric potential energy density \cite{forgacs2016}, which reads by
\begin{equation}
\mathcal{E}=\Sigma_{i=1}^2 |\textbf{D}\phi_i|^2+\frac{|\textbf{B}|^2}{2}+V(|\phi_i|), \label{1.1}
\end{equation}
where $\textbf{D}\phi_i=(\nabla-e_i\textbf{A})\phi_i,\ \textbf{B}=\nabla\times \textbf{A}$, with $e_i$ being the charge of the condensates $\phi_i$. The symmetric self-interaction potential of the two complex scalar fields is
\begin{equation}
V = \frac{\beta_1}{2}(|\phi_1|^2-1)^2+\frac{\beta_2}{2}|\phi_2|^4+\beta' |\phi_1|^2|\phi_2|^2-\alpha|\phi_2|^2.\label{1.2}
\end{equation}
This is a new type of solution from Witten's model \cite{witten1985}, superconducting strings with radially excited condensates. Radial excited states are common in non-topological soliton systems, such as Q-balls \cite{hartmann2017}, and topological soliton systems, such as magnetic monopoles \cite{mcleod2001}. \par

Different from the gauge cases, the ungauge (including pure scalar) cases are important, too. For the scalar field theory, it is always applied in high energy physics and if the symmetry is spontaneously broken globally, there must be a zero-mass boson, or Goldstone boson. And the remarkable feature of the global string \cite{vilenkin1994}, is the coupling with the Goldstone bosons, which will give rise to long range interactions. 
When the potential is symmetric, it is a simplification of the Standard model Higgs sector, and this is a 1VEV case. A non-symmetric case can be considered a toy model of the Higgs field of some symmetry breaking (electroweak or GUT) and the other scalar as some dark matter sector field

\par

There are few mathematical results for the TCGL theory. Alama, Bronsard and Mironescu considered the balanced case as the each solution goes to the same VEV with symmetric parameters in \cite{alama2006,alama2009,alama2012}. Cases for the general range of parameters has been done by Alama and Gao in \cite{qi2013} by a variational method, and stability in \cite{alama2014} by studying the spectrum of the second variational energy. Existence and uniqueness for constrained parameters and stability of symmetric vortices for the two-component $p$-Ginzburg-Landau model have been obtained in \cite{duan2020}.


In the present paper, we focus on the TCGL model because it is enough to expose the differences and to show the interesting vortices which are different from the ANO ones.
We provide a detailed existence proof, monotonicity, uniqueness, asymptotic description and quantization identity of the TCGL theory with extended scalar fields, established by P\'{e}ter Forg\'{a}cs and \'{A}rp\'{a}d Luk\'{a}cs in \cite{lukacsarpad} by using the shooting method and Schauder fixed point theorem. We consider the symmetry $U(1)\times U(1)$. The previous paper in physics focused on experiments with non-zero vacuum expectation values, which means that spontaneous symmetry breaking occurs in both fields and the system is transformed into a superconducting-superconducting state at the critical temperature, which is of significance for the study of multi-band superconductors. By choosing a different range of parameters, there is another situation: only one scalar field obtains nonzero vacuum expectation value. A non-zero VEV in the first field breaking the gauge symmetry lead to the production of vortices, while a vanishing VEV in the other field allows for it to condense in the core of the vortex of the first field. This situation corresponds to the superconducting-normal state.

Although part of our results (the 2VEV case) has already been considered in \cite{qi2013} by a variational method, we succeeded by shooting method which is more direct. The other part of things are different because the construction of the system corresponding to the 1VEV case, has no symmetry and the right hand side nonlinear term is changing sign, which can be regarded as a competitive system. And the asymmetry leads to the thin range of coefficients that we can choose. We do the existence and properties of solutions at the same time, by a constructive way. Since the system of equations are Euler equations and it is natural to use the integral function which indeed helps a lot (Greatly thanks to Mcleod and Wang \cite{mcleod2001}). Then, it can be constructed the set to do the Schauder fixed point theorem, where we establish the existence and monotonicity, with some other properties about the solutions.

In the following sections, we shall first state the model and our main results in Section 2 and study the existence, uniqueness, asymptotic properties and quantization of the 1VEV case in Section 3. Finally, we study the 2VEV case with existence and asymptotic properties in Section 4. (The uniqueness and quantization for the two cases is the same and we will claim them in Section 3).

\section{Model and main results}\label{s2}
We consider the complex scalar fields with $U(1)\times U(1)$ symmetric self-interaction potential
\begin{equation}
\label{1}
V(\phi_1,\phi_2) = \frac{\beta_1}{2}(|\phi_1|^2-1)^2+\frac{\beta_2}{2}|\phi_2|^4+\beta' |\phi_1|^2|\phi_2|^2-\alpha|\phi_2|^2,
\end{equation}
which includes $\phi_1$, $\phi_2$ representing the potentials that each scalar fields attain with 4 real parameters: $\beta_{1}$, $\beta_{2}$, $\beta^{'}$, $\alpha$. Throughout this paper, the following assumptions are made
\begin{equation}\beta_1>0, \quad\beta_2>0, \quad\beta'>-\sqrt{\beta_1 \beta_2},
\end{equation}
 to ensure $V$ is positive as $|\phi_{1}|^{2}$, $|\phi_{2}|^{2}\to \infty$.
Let us define Lagrangian
\begin{equation}
\label{2}
\mathcal{L} = \partial_{\mu}\Phi^{\dagger}\partial^{\mu}\Phi-V(\Phi^{\dagger}, \Phi),
\end{equation}
where $\Phi^{\dagger}=(\phi_1, \phi_2)^{\dagger}$, the potential $V$ is given by \eqref{1}. We choose the ansatz in \eqref{2} as
\begin{equation}
\phi_1=f(r)e^{in\theta}, \quad \phi_2=g(r)e^{im\theta},
\end{equation}
in which $n, m$ are non-negative integers describing the winding number or vortex number; $(r, \theta)$ are polar coordinates in the plane. Then the potential $V$ and the equations of motion for \eqref{2} can be reduced to the following
\begin{equation}
V_{ansatz} = \frac{\beta_1}{2}(|f|^2-1)^2+\frac{\beta_2}{2}|g|^4+\beta' |f|^2|g|^2-\alpha|g|^2,
\end{equation}
and
\begin{eqnarray}
f''+\frac{1}{r}f'-\frac{n^2}{r^2}f & = &f\left[\beta_{1}\left(f^{2}-1\right)+\beta^{'}g^{2}\right]   \label{f1},\\
g''+\frac{1}{r}g'-\frac{m^2}{r^2}g & = &g\left(\beta_{2}g^{2}-\alpha+\beta^{'}f^{2}\right) \label{f2},
\end{eqnarray}
respectively.

We see that it has three possible minima in \eqref{1}. Among them, $(|\phi_1|,|\phi_2|)=(1,0)$ and $(0,\sqrt{\alpha\mbox{/}\beta_2})$ are in one single scalar field of 1VEV. When the coefficients satisfy
\begin{equation}
\alpha>\beta',\quad {\beta'}^2<\beta_1\beta_2,
\end{equation}
and
\begin{equation}
(|\phi_1|,|\phi_2|)=\left(\sqrt{\frac{\beta_{1}\beta_{2}-\beta^{'}\alpha}{\beta_{1}\beta_{2}-\beta^{'^2}}},\sqrt{\frac{\beta_{1}(\alpha-\beta^{'})}{\beta_{1}\beta_{2}-\beta^{'^2}}} \right)\triangleq\left(A,B\right),
\end{equation}
will imply both fields acquire a non-zero VEV (2VEV). In the case of 2VEV, we can write the system into a symmetric form
\begin{eqnarray}
f''+\frac{1}{r}f'-\frac{n^2}{r^2}f & = &f\left[\beta_{1}\left(f^{2}-A\right)+\beta^{'}\left(g^{2}-B\right)\right]   \label{2f1},\\
g''+\frac{1}{r}g'-\frac{m^2}{r^2}g & = &g\left[\beta_{2}\left(g^{2}-B\right)+\beta^{'}\left(f^{2}-A\right)\right] \label{2f2},
\end{eqnarray}
and this is indeed the case Alama and Gao discuss. It is noted that the potential $V_{ansatz}$ is monotonous about distance between two-vortex configuration \cite{dantas2015,lukacsarpad}. \par
However, in the case of 1VEV, $(1,0)$ is of importance because that it is the global minimum. The case $(|\phi_1|,|\phi_2|)=(0,\sqrt{\alpha\mbox{/}\beta_2})$ can be obtained by changing parameters, so we just need to consider the previous case. Thus we would like to consider two different boundary conditions corresponding to the 1VEV case and 2VEV case for \eqref{f1}-\eqref{f2}.\par
\medskip
($i$)
For any $n\in \mathbb{N}^{+}$, $m=0$, $\beta_{1}, \beta_{2}>0,\ \beta^{'}>-\sqrt{\beta_1\beta_2}$ and there exists $\eta\in(0,1)$ such that
\begin{equation}
\frac{p_0}{\sqrt{\alpha-\beta'R^2\eta^{2\lambda}}}<\eta\label{c1}
\end{equation}
where $R$ is a positive constant, $p_0$ is the first zero of the classical Bessel equation. And the equation satisfy boundary conditions
{\setlength\arraycolsep{2pt} 
\begin{eqnarray}
&f(0)=0,\quad g'(0)=0, \label{bd3}\\
&f(\infty)=1,\quad g(\infty)=0,\label{bd4}
\end{eqnarray}}
where we expect a nontrivial real function $g$ with positive, finite initial value recorded by $b$ (need to be determined later), and we also assume
\begin{equation}
\beta_1\beta_2>\beta'^{2},\quad \beta^{\prime}> \alpha>0.\label{bc1}
\end{equation}

The above boundary conditions correspond to the 1VEV case. And the following boundary conditions describe the 2VEV case.\par
($ii$)
For any $m, n\in \mathbb{N}^{+}$, $\beta_{1}, \beta_{2}>0,\ \beta^{'}>-\sqrt{\beta_1\beta_2}$
{\setlength\arraycolsep{2pt} 
\begin{eqnarray}
&f(0)=g(0)=0,\label{bd1}  \\
&f(\infty)= A,\quad g(\infty)= B.\label{bd2}
\end{eqnarray}}
and we assume that
\begin{equation}
\beta_1\beta_2>\beta'^{2},\quad \alpha>\beta' \label{bc}
\end{equation}
clearly we have $\beta_1\beta_2>\beta'\alpha$ from \eqref{bc}.\\

\medskip

 In the following we state our main results for the boundary value problems \eqref{f1}-\eqref{f2} under the two different boundary value conditions \eqref{bd3}-\eqref{bd4} and \eqref{bd1}-\eqref{bd2}.

\begin{theorem}\label{2.1}
When $m=0$, for any positive integers $n>0$, given parameter condition \eqref{bc1} and \eqref{c1}, boundary value problems \eqref{f1}-\eqref{f2}, \eqref{bd3}-\eqref{bd4} admit a unique solution pair $(f, g)$, and each of them belongs to $C[0, \infty)$, satisfying
\begin{eqnarray*}
&f(r)&=D_0r^n+O(r^{n+2}), \quad g(r)= b_{0}+\frac{b_0(\beta_2 b_{0}^{2}-\alpha)}{4} r^{2} +O(r^{2\lambda+2}),\quad r\to 0^{+},\\
&f(r)&=1+O(r^{-2}), \quad g(r)=O(e^{-\sqrt{\beta'-\alpha}r}),\quad r\to\infty,\\
 &0\le& f(r)<1, \quad f'(r)>0,\quad  \forall r>0,\\
 &0<& g(r)\le b_0, \quad g'(r)<0,\quad  \forall r>0,
\end{eqnarray*}
where $b_0, D_0$ are positive constants, with $0<\lambda<n$. Also, the solution is quantized, i.e.
\begin{eqnarray}
&&\int_0^{\infty}rV_{ansatz}dr\nonumber\\
&=&\int_0^{\infty}r\left[\frac{\beta_1}{2}(|f|^2-1)^2+\frac{\beta_2}{2}|g|^4+\beta' |f|^2|g|^2-\alpha|g|^2\right]dr\nonumber\\
&=&\frac{ n^2}{2}.
\end{eqnarray}

\end{theorem}
\medskip
\begin{remark}\label{1.3}
We can regard this 1VEV case as a generalization of the classical GL theory. When $g\equiv0$, the system degenerate to the case that Chen, Elliott and Tang discussed in \cite{chen1994}. There are still some results like quantization effects \cite{brezis1994}.
\end{remark}

Noting that any solution of \eqref{f1}-\eqref{f2} are bounded when $r \to 0$, and it is a solution of the following equations
\begin{eqnarray}
f_{D}''+\frac{1}{r}f_{D}'-\frac{n^2}{r^2}f_{D} & = &f_{D}\big[\beta_{1}(f_{D}^{2}-1)+\beta^{'}g_{b}^{2}\big],\quad f_{D}\sim Dr^{n},\quad r\to 0^{+},   \label{f5}\\
g_{b}''+\frac{1}{r}g_{b}'-\frac{m^2}{r^2}g_{b} & = &g_{b}(\beta_{2}g_{b}^{2}-\alpha+\beta^{'}f_{D}^{2}),\quad g_{b}\sim b,\quad r\to 0^{+}, \label{f6}
\end{eqnarray}
with some positive constants $D$, $b$. In Theorem 2.1, the parameters we assume describe the $n^{th}$ derivative of $f$ and the initial value of $g$ at $0$, where $D$, $b$ are shooting parameters. To simplify, we denote $f_D$, $g_b$ by $f$, $g$.\par

The second result is the existence of unique solution in the 2VEV case, for any positive integers $m,\ n$.

\medskip
\begin{theorem}\label{2.2}
For any positive integers $m,\ n>0$, given the parameter condition \eqref{bc}, boundary value problems \eqref{f1}-\eqref{f2} and \eqref{bd1}-\eqref{bd2} admit a unique solution pair $(f, g)$, and each of them belongs to $C[0, \infty)$, satisfying
\begin{eqnarray*}
&f(r)&=D_0r^n+O(r^{n+2}),\quad g(r)=C_0r^m+O(r^{m+2}), \quad  r \to 0^+,\\
&f(r)&=A+O(r^{-2}),\quad g(r)=B+O(r^{-2}),\quad r\to\infty,\\
&0\le& f(r)<A, \quad f'(r)>0,\quad  \forall r>0,\\
&0\le& g(r)<B, \quad g'(r)>0, \quad \forall r>0,
\end{eqnarray*}
where $C_0$, $D_0$ are positive constants. Moreover, the solution to the potential is quantized
\begin{equation}
\int_0^{\infty}r(V_{ansatz}-V_{\infty})dr=\frac{n^2 A^2+m^2 B^2}{2},
\end{equation}
where $$V_{\infty}=\lim_{r\to \infty}V_{ansatz}(f(r),g(r))=V_{ansatz}(A,B)=\frac{\beta_1(\alpha-\beta')^2[\beta'(\beta'+\alpha)+(\beta'\alpha-\beta_1\beta_2)]}{2(\beta_1\beta_2-\beta'^2)^2}.$$.
\end{theorem}
\begin{remark}\label{2.2} When $m=n=0$, we see clearly that $f\equiv A$ and $g\equiv B$ are a solution to \eqref{f1}-\eqref{f2}. But we claim that it is the unique solution to the problem below from the uniqueness \cite{qi2013}. And it implies that no symmetry breaking has occurred.
\end{remark}

As we expect, the result is really an extension of the classical GL theory because of the symmetric of the system. And both of the components $f$ and $g$ are monotonically increasing to $A$ and $B$, their limiting value as $r\to\infty$, respectively.\par
Noting that any solution of \eqref{f1}-\eqref{f2} are bounded when $r \to 0$, and it is a solution to the following equations
\begin{eqnarray}
f_{D}''+\frac{1}{r}f_{D}'-\frac{n^2}{r^2}f_{D} & = &f_{D}\big[\beta_{1}(f_{D}^{2}-1)+\beta^{'}g_{C}^{2}\big],\quad f_{D}\sim Dr^{n},\quad r\to 0^{+},   \label{f3}\\
g_{C}''+\frac{1}{r}g_{C}'-\frac{m^2}{r^2}g_{C} & = &g_{C}(\beta_{2}g_{C}^{2}-\alpha+\beta^{'}f_{D}^{2}),\quad g_{C}\sim Cr^{m},\quad r\to 0^{+}, \label{f4}
\end{eqnarray}
with some positive constants $C$, $D$. In Theorem 2.2, the parameters we assume all describe the $n^{th}$ and $m^{th}$ derivative of $f$ and $g$ at $0$, respectively where $C$, $D$ are shooting parameters. To simplify, we denote $f_D$, $g_C$ by $f$, $g$.

It is well known that the existence of solutions is always transformed into the convergence of suitable sequences. In the following sections, we will use the shooting method for each function, to solve the non-emptiness, one ingeniously applied integral function and Schauder fixed point theorem, first used in \cite{mcleod2001}. In our case, there is a strong singularity near 0 different from the case in \cite{mcleod2001}, but we can overcome it by the special structure of the nonlinear terms.

\section{Existence of 1VEV case}\label{s3}
\hskip\parindent \baselineskip 0.2in
\renewcommand{\theequation}{3.\arabic{equation}}%
\setcounter{equation}{0}

In this section we prove Theorem 2.1. The proof is divided into two lemmas in which we choose to use integral function to get a natural sequence because the homogeneous part is Euler equation. We fix some $f$ in a suitable set and consider the equation for $g$ to get the solution $g_f$. Then we plug $g_f$ into the equation for $f$ and denote the solution by $\tilde{f}$. Finally, we prove that $\tilde{f}$ stay in the same set as $f$, and Theorem 2.1 follows from Schauder fixed point theorem.

\medskip
\begin{lemma}\label{3.1}
 Given any increasing function $f \in C[0, \infty)$ satisfying $f \le Rr^{\lambda}$, for $r \le 1$, together with $f(\infty)= 1$, we can find a unique, continuous and differentiable function $g$ satisfying \eqref{f2} and the conditions $g(0)=b_0$, $g(\infty)=0$. Moreover, $g$ is decreasing. Here, $\lambda$ is any fixed number with $0<\lambda<n$; $R$ is a positive constant, depending only on the given parameters $\alpha$, $\beta_{1}$, $\beta_{2}$, $\beta'$ in the problems.
\end{lemma}

\medskip

\noindent $\bf{Proof.}$ To make the proof clear, we divide it into six steps.\par
\noindent $\bf{Step.1}$ (Transform into the integral equation.)  Note that $m=0$, we convert \eqref{f2} into the integral equation by constant variation
\begin{equation}
g(r)=b+\hat{b}\ln r+\int_{0}^{r} s\left(\ln r-\ln s\right)g\left(\beta_2 g^{2}-\alpha+\beta^{\prime}f^{2}\right)ds,\quad r\in[0,\delta(b)).\label{intg1}
\end{equation}
By considering the corresponding boundary value condition, \eqref{intg1} becomes
\begin{equation}
g(r)=b+\int_{0}^{r} s\left(\ln r-\ln s\right)g\left(\beta_2 g^{2}-\alpha+\beta^{\prime}f^{2}\right)ds.\label{intg2}
\end{equation}
It is natural to get
\begin{equation}
g(r)= b+\frac{b(\beta_2 b^{2}-\alpha)}{4} r^{2} +O(r^{2\lambda+2}),\quad r\to 0,
\end{equation}
where we could conclude that the desired solution requires the range $b^2<\alpha\mbox{/}\beta_2$.\par
\noindent $\bf{Step.2}$ (Construction for the 'BAD' sets.)
Now, this boundary value problem only consists of one parameter $b$ which characterizes the initial value of $g$, and we interest in the case  $b>0$. In order to find a suitable solution, one constructs these two sets
\begin{align*}
& S_{1}\triangleq \left\{b>0|\ g' \ \text{becomes non-negative before }g\ \text{reaches }0\right\},\\
& S_{2}\triangleq \{b>0|\ g\ \text{crosses}\ 0\  \text{before}\ g'\ \text{becomes zero}\}.
\end{align*}
Theory of ordinary differential equation ensures that the solution has a continuous dependence on the initial value, so these two sets are both
 open. Also, the construction of $S_1$ and $S_2$ guarantees that they are disjointed. Then we claim that both of them are non-empty.\par

\noindent $\bf{Step.3}$ ($S_1$ is nonempty.)
We derive \eqref{intg2} and obtain
\begin{equation}
g'(r)=b+\int_{0}^{r} \frac{s}{r} g\left(\beta_2 g^{2}-\alpha+\beta^{\prime}f^{2}\right)ds.
\end{equation}
If we choose $b$ sufficiently large, $g'(r)$ stays positive for any $r\ge 0$.
Hence $S_{1}$ is non-empty.\par

\noindent $\bf{Step.4}$ ($S_2$ is nonempty.)
Let us define $w=\frac{g}{b}$, then $w(0)=1$ and $w$ satisfies
  \begin{equation}
  w''+\frac{1}{r}w'=w\left(\beta_2b^2w^2-\alpha+\beta'f^2\right).\label{w}
  \end{equation}
 When $b\to 0$, equation \eqref{w} converts to
 \begin{equation}
 w''+\frac{1}{r}w'=w\left(-\alpha+\beta'f^2\right).\label{w1}
 \end{equation}\par
 We denote $p_0$ to be first zero of the classical Bessel equation
 \begin{gather*}
 p''+\frac{1}{r}p'+p=0,\\
 p(0)=1, \quad p'(0)=0.
 \end{gather*}

  We want to claim $S_2$ contains small $b$ by analysing equation \eqref{w1}, which is similar to the oscillating zero order Bessel equation if we analysis it in a proper interval. The difficulty is that term $(-\alpha+\beta'f^2)$ changes sign in $(0,\infty)$ and this is the reason why we make the assumption for the coefficients. \par
  We can assume that $w>0$ in $[0,1]$ and if not, the proof is finished.  We note that when $r\in[0,1]$, $f\le Rr^{\lambda}$, and $\beta'>\alpha$, there exists $\eta\in(0,1)$ such that
  \begin{equation}
  -\alpha+\beta'R^2\eta^{2\lambda}<-\left(\frac{p_0}{\eta}\right)^2,
  \end{equation}
  and then
  \begin{equation}
  w''+\frac{1}{r}w'=w\left(-\alpha+\beta'f^2\right)\le w\left(-\alpha+\beta'R^2\eta^{2\lambda}\right).
  \end{equation}
  We consider the comparing equation
  \begin{gather}
  \hat{w}''+\frac{1}{r}\hat{w}'=\hat{w}\left(-\alpha+\beta'R^2\eta^{2\lambda}\right),\\
  \hat{w}(0)=1,\quad \hat{w}'(0)=0,
  \end{gather}
  and it is an oscillating zero order Bessel equation with the first zero
  \begin{equation}
  \hat{w}_0=\frac{p_0}{\sqrt{\alpha-\beta'R^2\eta^{2\lambda}}}\in(0,\eta).
  \end{equation}
  Then we know that $w$ will cross zero as $b$ small enough and $S_2$ is nonempty.

\noindent $\bf{Step.5}$ (The existence of $g$.)

However, $\mathbb{R}^{+}=(0,\infty)$ cannot consist of two open, nonempty, disjointed sets, there must exist $b_0$ neither in $S_1$ nor $S_2$, which satisfy
\begin{equation}
0< b_0<\sqrt{\frac{\alpha}{\beta_2}} < \sqrt{\frac{\alpha+|\beta'|}{\beta_2}}.
\end{equation}
Moreover, the closed set where $b_0$ is in, has these properties
\begin{eqnarray}
 g'(r)< 0,\quad 0< g(r)\le b_0,\quad for\ all \ r\ge 0. \label{b0}
\end{eqnarray}
So, the solution exists for any $r\ge0$ and satisfies the boundary condition $\lim_{r\to\infty}g(r)=0$. Otherwise,
\begin{equation}
\frac{1}{r}\left(rg_r\right)_r\neq 0,\quad for\ all \ r\ge 0.
\end{equation}
It is not true from \eqref{b0}.\par

\noindent $\bf{Step.6}$ (Uniqueness.)
 If there are $g_1$ and $g_2$ satisfying the boundary conditions
\begin{eqnarray}
&g_1(0)=b_1,\quad g_2(0)=b_2,\\
&g_1'(0)=g_2'(0)=0,\quad g_1(\infty)=g_2(\infty)=0,
\end{eqnarray}
where $b_1\not= b_2$.
Note that for any $i=1,2$, we obtain the estimate
\begin{equation}
g_i=O(e^{-\sqrt{\beta'-\alpha}r}),\quad r\to\infty.
\end{equation}
Let us rearrange the equation
\begin{equation}
\frac{(rg_i')'}{g_i}=r(\beta_2g_i^{2}-\alpha+\beta'f^2).\label{reg}
\end{equation}
Subtracting the two equations, multiplying $g_2^2-g_1^2$ on the both sides, then integrating from $0$ to infinity which is the idea from Brezis and Oswald in \cite{brezis1986}, we obtain
\begin{eqnarray}
\int_0^{\infty}\left(\frac{(rg_1')'}{g_1}-\frac{(rg_2')'}{g_1}\right)\left(g_2^2-g_1^2\right)dr&=&\int_0^{\infty}r\left(\left|g_1'-\frac{g_2'g_1}{g_2}\right|^2+\left|g_2'-\frac{g_1'g_2}{g_1}\right|^2\right)dr \nonumber\\
&=&-\int_0^{\infty}r\beta_2\left(g_1^2-g_2^2\right)^2dr\le 0.
\end{eqnarray}
It means that we have
\begin{eqnarray}
g_1'g_2-g_2'g_1=g_2^2\left(\frac{g_1}{g_2}\right)'=0,
\end{eqnarray}
i.e. there are two cases: $g_2\equiv 0$ or $g_1=kg_2$, where $k$ is a constant. We get a contradiction because when we plug it into \eqref{reg} we get that $k=\pm 1$. Recalling that $g_i$ are both positive and we prove the uniqueness.

Clearly, $b_0$ has a uniform upper bound
\begin{equation}
0<g(r)\le b_0<\sqrt{\frac{\alpha}{\beta_2}},\quad for\ all \ r\ge 0,
\end{equation}
which is useful in the proof of the next lemma.\qed

\medskip

\begin{lemma}\label{3.2}
For the $g$, actually $g(f)$ obtained in Lemma 3.1, we can find a unique $\tilde{f}\in C[0, \infty)$, satisfying
\begin{equation}
\tilde{f}''+\frac{1}{r}\tilde{f}'-\frac{n^2}{r^2}\tilde{f}=\tilde{f}\left[\beta_1(\tilde{f}^2-1)+\beta'g^2\right],   \label{tilf1}
\end{equation}
and
\begin{equation}
\tilde{f}(0)=0, \quad \tilde{f}(\infty)=1,
\quad 0\leq \tilde{f}(r)<1, \quad \tilde{f}'(r)>0, \quad 0\le r\le \infty.
\end{equation}
Moreover, for the $R$ in Lemma 3.1, we have $\tilde{f}\le R^{*}(R)r^{n},\ r\le 1$.
\end{lemma}

\medskip

\noindent $\bf{Proof.}$To make the proof clear, we divide it into six steps.\par
\noindent $\bf{Step.1}$
 The fundamental solutions of \eqref{tilf1} are $r^n, r^{-n}$. We can convert it formally into the integral
equation from the local existence of ordinary differential equation
\begin{equation}
\tilde{f}(r)=Dr^{n}+\frac{1}{2n}\int_{0}^{r}
\left(\frac{r^{n}}{s^{n-1}}-\frac{s^{n+1}}{r^{n}}\right)\tilde{f}\left[\beta_1(\tilde{f}^2-1)+\beta'g^{2}\right]ds, \quad r\in[0,\delta(D)).\label{tilfn}
\end{equation}
Equation \eqref{tilfn} can be solved by iteration, at least for $r$ sufficiently small, where $\tilde{f}$ has a natural estimate
\begin{equation}
\tilde{f}(r)=Dr^{n}+O(r^{n+2}).\label{tilfn0}
\end{equation}
We are interested in the case $D>0$. Let us define the following sets
{\setlength\arraycolsep{2pt}
\begin{align*}
& S_{3}\triangleq \left\{D>0\big|\tilde{f}'\ \text{reaches} \ 0\ \text{before}\ \tilde{f}\ \text{reaches}\ 1\right\},\\
& S_{4}\triangleq \left\{D>0\big|\tilde{f}\ \text{crosses}\ 1\ \text{before}\ \tilde{f}'\ \text{reaches}\ 0\right\}.
\end{align*}}
The solution depends only on the parameter $D$ continuously, so they are both open sets. From the constructions of $S_3$ and $S_4$, we know that $S_{3}\cap S_{4}=\varnothing$.\par

\noindent $\bf{Step.2}$
First, we claim that $S_3$ is  non-empty. Noting that when $D=0$, we get $\tilde{f}\equiv 0$. Then, by the continuous dependence on parameter, as $D$ sufficiently small, $\tilde{f}_D$ is small, too. Let $\tilde{f}=Du$ in \eqref{tilf1}, then $u$ satisfies
\begin{equation}
u''+\frac{1}{r}u'-\frac{n^2}{r^2}u = u\left[\beta_1(D^2u^2-1)+\beta'g^2\right]
. \label{un1}
\end{equation}
Let $D\to 0$, by following the continuity, the term $D^2u^2$ is a higher order infinitesimal. Equation \eqref{un1} becomes
\begin{equation}
u''+\frac{1}{r}u'-\frac{n^2}{r^2}u = u(\beta'g^2-\beta_1)\label{u}.
\end{equation}
We interest in $u>0$, otherwise, $S_3$ is non-empty obviously. Rewrite equation \eqref{u} as
{\setlength\arraycolsep{2pt}
\begin{eqnarray}
\begin{split}
u''+\frac{1}{r}u'+\left(\beta_1-\beta'b_0^2-\frac{n^2}{r^2}\right)u&=u\beta'(g^2-b_0^2)<0, \label{un2} \\
u\sim r^{n},\quad r\to 0&.
\end{split}
\end{eqnarray}}
We consider the comparison equation of \eqref{un2}
{\setlength\arraycolsep{2pt}
\begin{eqnarray}
\begin{split}
v''+\frac{1}{r}v&'+\left(\beta_1-\beta'b_0^2-\frac{n^2}{r^2}\right)v=0, \label{un3} \\
&v\sim r^n,\quad r\to 0.
\end{split}
\end{eqnarray}}
Assume that $v\ge 0,\ r\in [0, R_{b_0})$. We can conclude that $v>u$ in $(0, \delta)$ when $\delta$ is sufficiently small. We claim that $v>u$, $r\in[0,R_{b_0})$. If not, there exists $r_1\in (0,R_{b_0})$, such that $v(r_1)=u(r_1)$, one constructs
\begin{equation}
c=\sup\left\{\hat{c}>0|\ \hat{c}v\le u,\ r\in (0, r_1)\right\},
\end{equation}
which implies that there exists $r_2\in (0,r_1)$, satisfying
\begin{equation}
cv=u,\quad cv'=u',\quad cv''\le u'',\quad r=r_2;
\end{equation}
but it is impossible, because
{\setlength\arraycolsep{2pt}
\begin{eqnarray}
0&=&\left.c\left[v''+\frac{1}{r}v'+\left(\beta_1-\beta'b_0^2-\frac{n^2}{r^2}\right)v\right]\right|_{r=r_2}\nonumber\\
&\le&\left.u''+\frac{1}{r}u'+\left(\beta_1-\beta'b_0^2-\frac{n^2}{r^2}\right)u\right|_{r=r_2}\nonumber\\
&=&u(r_2)\beta'\left(g^2(r_2)-b_0^2\right)<0.
\end{eqnarray}}
A contradiction! So $v>u>0,\ r\in[0, R_1)$, $\forall\ R_1 \le R_{b_0}$. Moreover, recalling that $\beta_1-\beta'b_0^2>0$ and it is the oscillating Bessel equation. Thus there must exist $r_3$ such that $u'(r_3)=0$, which means $S_3$ is non-empty.

\medskip

\begin{proposition} If there exists $r'$ such that it is the first one satisfying $\tilde{f}(r')=1$, then $\forall\ r>r'$, we conclude that $\tilde{f}'(r)>0$.
\end{proposition}
\medskip

\noindent $\bf{Proof.}$ In this case, it is similar to the standard Ginzburg-Landau equation, and we recall the proof of lemma 2.6 in \cite{chen1994}. If $\tilde{f'}(r')=0$, then
\begin{equation}
\tilde{f}''(r')=\lim_{h\to 0}\frac{\tilde{f}'(r')h-(\tilde{f}(r'-h)-\tilde{f}(r'-2h))}{h^{2}}\le 0,
\end{equation}
but we get $\tilde{f}''(r')=\tilde{f}(\beta'g^{2}+{n^{2}}\mbox{/}{r'^{2}})>0$, which is a contradiction. So $\tilde{f}'(r')>0$. If $\forall\ r>r'$, $\tilde{f}'(r)\le 0$, then there exists $r''$ such that $\tilde{f}'(r'')=0$ and $\tilde{f}''(r'')\le 0$. Let $r=r''$ in \eqref{tilf1} and it will not happen.\qed

\medskip
\par

\noindent $\bf{Step.3}$
Next we will prove that $S_4$ is non-empty. From \eqref{tilfn}, we recall that
\begin{eqnarray}
\tilde{f}(r)&=&Dr^{n}+\frac{1}{2n}\int_{0}^{r}\left(\frac{r^{n}}{s^{n-1}}-\frac{s^{n+1}}{r^{n}}\right)D^{-\frac{2}{n}}\tilde{f}\left[\beta_1(\tilde{f}^2-1)+\beta'g^{2}\right]d s. \nonumber\\
&\ge&Dr^{n}+\frac{1}{2n}\int_{0}^{r}\left(\frac{r^{n}}{s^{n-1}}-\frac{s^{n+1}}{r^{n}}\right)D^{-\frac{2}{n}}\tilde{f}\left[\beta_1(\tilde{f}^2-1)\right].
\label{tilfnt}
\end{eqnarray}
If we choose the initial iterative function $\tilde{f}_0=t^n$, then we can extend it in any compact interval (it means that $\tilde{f}(t)$ cannot blowup at any finite points). We have
\begin{align}
\tilde{f}(t)&=t^n+\frac{D^{-\frac{2}{n}}}{2n}\int_0^{t}\tau t^n(1-\frac{\tau^{2n}}{t^{2n}})\frac{\tilde{f}}{\tau^{n}}\left[\beta_1(\tilde{f}^2-1)+\beta'g^{2}\right]d \tau,\label{int}
\end{align}
The only possible singularity occurs at $t=0$, but it is a removable singularity. In fact,
\begin{equation}
\lim_{t\to 0}\frac{\tilde{f}}{t^n}=1.
\end{equation}
I. e. for any fixed $\epsilon>0$, there exists $\delta>0$, such that for any $t\in [0,\delta]$, we get
\begin{equation}
0<\left|\frac{\tilde{f}}{t^n}\right|<1+\epsilon.
\end{equation}
Now, Let us estimate the convergence of the integration as $t\in [0,\delta]$ to clarify if in any compact interval, we can extend the iteration, for any fixed $D$,
\begin{align}
&\left|\frac{D^{-\frac{2}{n}}}{2n}\int_0^{t}\tau t^n(1-\frac{\tau^{2n}}{t^{2n}})\frac{\tilde{f}}{\tau^{n}}\beta'g^{2}d\tau\right|\nonumber\\
\le&\frac{D^{-\frac{2}{n}}}{2n}\int_0^{t}\left|\tau t^n\frac{\tilde{f}}{\tau^{n}}\beta'g^{2}\right|d\tau\nonumber\\
\le&\frac{D^{-\frac{2}{n}}\beta'b_0^2\delta^n}{2n}\int_0^{\delta}(1+\epsilon)\tau d\tau \rightarrow 0,\quad \text{as}\ \delta\to 0,
\end{align}
where we could conclude that \eqref{int} can be extended in any finite interval. Then, equation \eqref{tilfnt} is equivalent to the following equation as $D\to\infty$, in any compact subinterval of $(0,\infty)$
\begin{equation}
\tilde{f}_{tt}+\frac{1}{t}\tilde{f}_{t}-\frac{n^2}{t^2}\tilde{f}=D^{-\frac{2}{n}}\beta_{1}\tilde{f}(\tilde{f}^{2}-1).
\end{equation}
Making a formal power series expansion for $\tilde{f}$ when $t\in(0,\frac{\delta}{2})$, we get
\begin{align}
\tilde{f}&\sim t^n-\frac{D^{-\frac{2}{n}}\beta_1}{4n+4}t^{n+2}+o\left(t^{n+2}\right) \nonumber\\
           &\sim |D|\left(r^{n}-\frac{\beta_{1}}{4n+4} r^{n+2}\right)+o\left(r^{n+2}\right).
\end{align}
So there exists some positive constant $N$, satisfying $|D|^{\frac{2}{n}}>N\beta^{\prime}b^{2}$, where we can find a sufficiently large $D$ because of the continuity on the parameter, such that $\tilde{f}$ crosses $1$ when $r$ is sufficiently small and then $S_4$ is non-empty.
\par

\noindent $\bf{Step.4}$
If $D_0\notin S_3\cup S_4$, then the solution related to $D_0$ satisfies these properties
\begin{equation}
 0\le \tilde{f}<1,\quad \tilde{f}'>0.\quad for\ all\ r\ge 0.
\end{equation}
$\tilde{f} $ is increasing with an upper bound, so we can denote that $\lim_{r \rightarrow \infty}\tilde{f}(r)=\tilde{f}_{\infty}$, then $\tilde{f}_{\infty}\in (0, 1]$. Let $r\to \infty$ in the both sides of the equation \eqref{tilfn}, we obtain
\begin{equation}
\frac{1}{r}\left(r\tilde{f}_{r}\right)_{r}=\beta_1\tilde{f}_{\infty}(\tilde{f}^{2}_{\infty}-1),
\end{equation}
and conclude that $\tilde{f}_{\infty}=1$.
\par

\noindent $\bf{Step.5}$
It is required to prove the uniqueness of $\tilde{f}$ to form the map. Assume that there exists $\tilde{f}_1$ and $\tilde{f}_2$  satisfying
$\tilde{f}_1(0)=\tilde{f}_2(0)=0,\  \tilde{f}_1(\infty)=\tilde{f}_2(\infty)=1$. Let $q_i=r^{\frac{1}{2}}\tilde{f}_i$, $i=1,2$ in \eqref{tilfn}. It is transformed into
\begin{equation}
r^2q_i''=q_i\left[r\beta_1(q_i^2-r)+\beta'g^2r^2+n^2-\frac{1}{4}\right].\label{qn}
\end{equation}
Denote $Q=q_1-q_2$, which gives us
\begin{equation}
r^2Q''=Q\left[r\beta_1\left(q_1^2+q_1q_2+q_2^2-r\right)+\beta'g^2r^2+n^2-\frac{1}{4}\right].\label{Qn}
\end{equation}
 We know that $Q$ oscillates more slowly than $q_1$ from the Sturm comparison theorem, but $q_1$ only equals to 0 when $r=0$, so $Q$ can
not equal to 0 at any finite point. Without loss of generality, assume $Q>0$ initially.
Comparing \eqref{qn} and \eqref{Qn}, we know that ${Q}\mbox{/}{q_1}$ is increasing monotonically, then it should have a strictly positive limit, which is a contradiction and one gives the uniqueness of $\tilde{f}$.
\par

\noindent $\bf{Step.6}$
Lemma 3.2 is proved apart from the estimate of $\tilde{f}(r)$, $r\le 1$. From \eqref{tilfn}, we get
\begin{align}
\frac{\tilde{f}}{r^{n}}&=D_0+\frac{1}{2n}\int_{0}^{r}\left(\frac{1}{s^{n-1}}-\frac{s^{n+1}}{r^{2n}}\right)\tilde{f}\left[\beta_1(\tilde{f}^{2}-1)+\beta'g^{2}\right]ds, \nonumber\\
&=D_0+\frac{1}{2n}\int_{0}^{r}\left(s-\frac{s^{2n+1}}{r^{2n}}\right)\frac{\tilde{f}}{s^{n}}\left[\beta_1(\tilde{f}^{2}-1)+\beta'g^{2}\right]ds.\label{tilfr2}
\end{align}
and recalling that $\lim_{r\to 0} r^{-n}\tilde{f}=D_0$, thus we could know that there exists a uniform $M>0$ such that $r^{-n}\tilde{f}\le M$, $r\le 1$ which depends on $(N,R,\alpha, \beta',\beta_1,\beta_2)$. Also, the term $|\beta_1(\tilde{f}^{2}-1)+\beta'g^{2}|$ has a bound in $[0,1]$ because of the continuity in compact interval $[0,1]$, and one denotes it as $\hat{M}(R)$
\begin{eqnarray}
|\beta_1(\tilde{f}^{2}-1)+\beta'g^{2}|&\le& \beta_1\tilde{f}^{2}+\beta_1+|\beta'|g^{2}\nonumber\\
&\le&\beta_1\left(R^2r^{2\lambda}+1\right)+\frac{|\beta'|\alpha}{\beta_2}\triangleq \hat{M}(R).
\end{eqnarray}
 Then we have
\begin{eqnarray}
\left|\tilde{f}(r)-D_0r^{n}\right|&\le&\frac{1}{2n}\int_{0}^{r}{r^{n}s^{1-n}\tilde{f}\left|\beta_1(\tilde{f}^2-1)+\beta'g^2\right|}ds \nonumber\\
&\le&\frac{M\hat{M}(R)}{2n}r^{n}\int_{0}^{r} s ds\nonumber\\
&=&\frac{M\hat{M}(R)}{4n}r^{n+2},\quad 0\le r\le 1,
\end{eqnarray}
which means
\begin{eqnarray}
\left|\frac{\tilde{f}}{r^n}\right|&\le& D_0+\frac{1}{2n}\int_0^r M\hat{M}(R)sds\nonumber\\
&=&D_0+\frac{M\hat{M}(R)}{4n}r^2 \nonumber\\
&\le& D_0+\frac{M\hat{M}(R)}{4n}\triangleq R^*,
\end{eqnarray}
and $|r^{-n}\tilde{f}|\le R^{*}$ where $R^*$ depends on $(N,R,\alpha, \beta',\beta_1,\beta_2)$. The proof of Lemma 3.2 is carried out.\qed

\begin{remark}\label{3.1}
For any $f$ satisfying $f\le Rr^{\lambda},\ r\le1$, we gain that $\tilde{f}\le R^*r^n,\ r\le 1$. It is known that
\begin{equation}
\tilde{f}\le R^*r^n \le R r^{\lambda},\quad r\le \left(\frac{R}{R^*}\right)^{\frac{1}{n-\lambda}}\triangleq r_0.
\end{equation}
\end{remark}
\medskip

\noindent $\bf{Proof\ of\ Theorem\ 2.1.}$
Define the $Banach$ space $\mathscr{B}$
\begin{equation}
\mathscr{B}=\left\{f\in C[0, \infty)\ \big| \ \exists M>0,\ r^{-\lambda}(1+r^{\lambda})f\le M\right\},\label{banach}
\end{equation}
with the norm
\begin{equation}
\|f\|_{\mathscr{B}}=\sup_{r\in [0, \infty)}\left\{\big|r^{-\lambda}(1+r^{\lambda})f\big|\right\}.\label{norm}
\end{equation}
One chooses a non-empty bounded closed convex set $S$ as
\begin{equation}
S=\big\{f\in B:\big|r^{-n}f\big|\le R^{*}, \quad r\le r_0; \quad f(\infty)=1,\quad 0\le f< 1, \quad f\ \text{is increasing}\big\}.\label{set}
\end{equation}
Evidently, $S$ is non-empty, bounded, closed and convex. We need to prove the following\par
($i$) The mapping $\Phi$ maps $S$ into itself, and $\Phi:f\mapsto \tilde{f}$.\par
($ii$) $\Phi$ is continue.\par
($iii$) $\Phi$ is precompact.\par
Then the Schauder fixed point theorem ensures that $\Phi$ has at least one fixed point, and the existence of the system is proved.

Part($i$) is guaranteed by Lemma 2.1 and Lemma 2.2, and part($ii$) follows from continuous dependence of the solution on the parameters. Although the continuous dependence based on the maximum norm in continuous function space, we can prove that the norm $\|\cdot \|_{\mathscr{B}}$ we define is equivalent to it, i.e.
\begin{equation}
\|f\|_{L^{\infty}}\le \|f\|_{\mathscr{B}}=\sup_{r\in [0, \infty)}\left\{\big|r^{-\lambda}(1+r^{\lambda})f\big|\right\}\le 2\|f\|_{L^{\infty}}.
\end{equation}
Thus part($ii$) is done. To prove part($iii$), given any sequence $\{f_k\}\in S$, which is bounded in the norm of $\mathscr{B}$, then we see that $\tilde{f}_k'$ is bounded in any finite interval $(0,r)$. From the differential mean value theorem, one knows that $\tilde{f}_{k}$ is equicontinuous and there exists a subsequence of $\{\tilde{f}_k\}$, converging uniformly to some $\tilde{f}$ in any compact subinterval of $(0, \infty)$ by the Ascoli-Arzela theorem. The difficulties occur at $r\to 0$ and $r\to \infty$.

When $r\to 0$, given $\epsilon >0$, there exists $\hat{r}$ sufficiently small, such that
\begin{equation}
\sup_{r\le r_0}|r^{-\lambda}(\tilde{f_k}-\tilde{f})|=\sup_{r\le r_0}|r^{-n}(\tilde{f_k}-\tilde{f})r^{n-\lambda}|\le \hat{r}^{n-\lambda}
R^{*}<\epsilon.
\end{equation}

When $r\to \infty$, recall \eqref{tilf11} to consider the asymptotic property for $\tilde{f}_k$
\begin{equation}
\tilde{f}_k''+\frac{1}{r}\tilde{f}_k'-\frac{n^2}{r^2}\tilde{f}_k=\tilde{f}_k\left[\beta_1(\tilde{f}_k^2-1)+\beta'g^2\right].\label{tilf4}
\end{equation}
To obtain the fact we would like to get the estimate for $\tilde{f}_k$ as $r\to \infty$, let
\begin{equation}
e(r)=\beta_1 \tilde{f}_k(\tilde{f}_k+1),
\end{equation}
 then consider comparison function
\begin{equation}
F=\mu r^{-\sigma},\quad r>0,\quad \mu>0.\label{F}
\end{equation}
We get
\begin{equation}
F''+\frac{1}{r}F'-\frac{n^{2}}{r^{2}}F=\frac{\sigma^{2}-n^{2}}{r^{2}}F.\label{F1}
\end{equation}
and by combining with \eqref{tilf4}, one gets
\begin{eqnarray}
(F+\tilde{f}_k-1)''+\frac{1}{r}(F+\tilde{f}_k-1)'&=&\left[\frac{\sigma^{2}-n^{2}}{r^{2}}+\frac{n^{2}}{\mu}r^{\sigma -2}-e(r)\right]F+\beta' \tilde{f}_kg^{2} \nonumber \\
&&+(e(r)+\frac{n^{2}}{r^{2}})(F+\tilde{f}_k-1).
\end{eqnarray}
Set $\sigma=2$, and note that $g$ converges to $0$ exponentially
\begin{equation}
g(r)=O(e^{-\sqrt{\beta'-\alpha}r}), \quad r\to \infty.\label{estg}
\end{equation}
One finds that if we choose $\mu$ large enough, there exists $\bar{r}>1$, such that
\begin{equation}
\frac{4-n^{2}}{r^{2}}+\frac{n^{2}}{\mu}-e(r)+\frac{\beta'}{r^{3}}<0,
\end{equation}
and
\begin{equation}
(F+\tilde{f}_k-1)(\bar{r})>0.
\end{equation}
Thus we conclude
\begin{equation}
(F+\tilde{f}_k-A)''+\frac{1}{r}(F+\tilde{f}_k-1)'\le (e(r)+\frac{n^{2}}{r^{2}})(F+\tilde{f}_k-1).
\end{equation}
From maximum principle, we get
\begin{equation}(F+\tilde{f}_k-1)(r)\ge 0,\quad r\ge \bar{r};
\end{equation}
i.e.
\begin{equation}
-\mu r^{-2}\le \tilde{f}_k(r)-1 \le 0.
\end{equation}
For the upper bound, recalling \eqref{estg}, we can choose $\hat{\mu}$ small, such that
\begin{equation}
(F+\tilde{f}_k-A)(\bar{r})<0
\end{equation}
and
\begin{equation}
\left[\frac{4-n^{2}}{r^{2}}+\frac{n^{2}}{\hat{\mu}}-e(r)\right]F+\beta' \tilde{f}_kg^{2}\ge 0,\quad r>\bar{r}.
\end{equation}
By using Maximum principle again, we get
\begin{equation}
\tilde{f}_k\le 1-\hat{\mu} r^{-2}.
\end{equation}
We derive that
\begin{equation}
\tilde{f}_k=1+O(r^{-2}). \quad as\ r\to\infty. \label{estf}
\end{equation}
From all of these, we finish the uniform convergence as $r\ge 0$. This is indeed what we need. However, Schauder fixed point theorem can only guarantee the existence of the solution, but not the uniqueness of the system.
\par
\medskip
Uniqueness of the system follows from Alama and Gao in \cite{qi2013} which is the idea of Brezis and Oswald \cite{brezis1986}, and we briefly show here. If we have $(f_1, g_1)$ and $(f_2,g_2)$ both are the solutions to the system and corresponding boundary value conditions, we rearrange the equations and obtain
\begin{eqnarray}
\frac{(rf_i')'}{f_i}-\frac{n^2}{r^2}&=&r\left[\beta_1(f_i^2-1)+\beta'g_i^2\right], \quad i=1,2\label{r1}\\
\frac{(rg_i')'}{g_i}-\frac{m^2}{r^2}&=&r\left(\beta_2g_i^2-\alpha+\beta'f_i^2\right),\quad i=1,2.\label{r2}
\end{eqnarray}
For \eqref{r1}-\eqref{r2}, multiply $f_2^2-f_1^2$, $g_2^2-g_1^2$, respectively, integrate from $0$ to infinity, and plus them, we get
\begin{align}
&\int_0^{\infty}r\left(\left|f_1'-\frac{f_2'f_1}{f_2}\right|^2+\left|f_2'-\frac{f_1'f_2}{f_1}\right|^2+\left|g_1'-\frac{g_2'g_1}{g_2}\right|^2+\left|g_2'-\frac{g_1'g_2}{g_1}\right|^2\right)dr,\nonumber\\
=&-\int_0^{\infty}r\left[\beta_1(f_1^2-f_2^2)^2+2\beta'(f_1^2-f_2^2)(g_1^2-g_2^2)+\beta_2(g_1^2-g_2^2)^2\right]dr\le 0.
\end{align}
Then we proved the uniqueness.\par
\medskip
Now, let us talk about the quantization of the radial solution for general integer $n$, $m$ and vacuum expectation value.
Multiplying \eqref{f1}-\eqref{f2} by $r^2 f'$ and $r^2 g' $, respectively and integrating them over $(0,R)$, one leads to
\begin{equation}
\frac{R^2}{2}\left(f'^{2}(R)+g'^{2}(R)\right)-\left(\frac{n^2}{2}f^{2}(R)+\frac{m^2}{2}g^{2}(R)\right)=\frac {R^2}{2}\left(V(R)-c\right)-\int_0^R r(V-c)dr,
\end{equation}
then we get
\begin{equation}
\int_0^R r(V_{ansatz}-V_{\infty})dr=\frac{1}{2}\left(n^2f^{2}(R)+m^2g^{2}(R)\right)-\frac{R^2}{2}\left(f'^{2}(R)+g'^{2}(R)\right)+\frac {R^2}{2}\left(V_{ansatz}-V_{\infty}\right),\label{q}
\end{equation}
where $V_{\infty}=\lim_{r\rightarrow \infty}V_{ansatz}(f(r),g(r))$. Clearly, the relative derivative terms have the estimate
\begin{equation}
\lim_{R\to\infty}R^2 \left(f'(R)+g'(R)\right)=0,
\end{equation}
from \eqref{estg} and \eqref{estf}. Next, we know
\begin{eqnarray}
\lim_{R\rightarrow\infty} \left(V_{ansatz}(R)-V_{\infty}\right)&=&\frac{\beta_1}{2}\left(f^2-A^2\right)^2+\frac{\beta_2}{2}\left(g^2-B^2\right)^2+\beta'\left(f^2-A^2\right)\left(g^2-B^2\right),\nonumber\\
&=&O\left(\frac{1}{R^4}\right).
\end{eqnarray}
So, this quantized identity is acquired by
\begin{equation}
\int_0^R r(V_{ansatz}-V_{\infty})dr=\frac{1}{2}\left(n^2A^{2}+m^2B^{2}\right),\quad R\to\infty,
\end{equation}
which implies
\begin{equation}
\lim_{R\to\infty}\int_0^{2\pi}\int_0^R r(V_{ansatz}-V_{\infty})drd\theta=\pi\left(n^2A^{2}+m^2B^{2}\right),
\end{equation}
for any positive integers $n$, $m$ and we prove Theorem 2.1.\qed

\section{Existence of the 2VEV case}\label{s4}

\setcounter{equation}{0}\setcounter{remark}{0}
\setcounter{theorem}{0}\setcounter{remark}{0}
\setcounter{lemma}{0}\setcounter{remark}{0}

We prove Theorem 2.2 in this section. The only difference between the two cases is the shooting parameters we choose. Noting that the structure of the system, which is symmetric, the proof will be simplified. Theorem 2.2 is also proved by Schauder fixed point theorem.

\medskip

\begin{lemma}\label{4.1}
Given any function $f \in C[0, \infty)$ such that $f \le Rr^{\lambda}$, for $r \le 1$, $f$ is increasing in $\mathbb{R}^+$ and $f(\infty)= A$, we can find a unique, continuous and differentiable function $g$ satisfying \eqref{f2} and $g(0)=0,\ g(\infty)=B$. Moreover, $g$ is increasing. Here, $\lambda$ is any fixed number with $0<\lambda<n$; $R$ is a large, positive constant, depending only on the given parameters $(\alpha, \beta_{1}, \beta_{2}, \beta')$ in the equations.
\end{lemma}

\medskip

\noindent $\bf{Proof.}$ We divide the proof into six steps.
\par

\noindent $\bf{Step.1}$
We note the fact that $r^{m}, r^{-m}$ are fundamental solutions of homogeneous equation corresponding to \eqref{f2} , we can formally write it into the integral equation when $r$ is sufficiently small
\begin{equation}
g(r)=Cr^{m}+\frac{1}{2m}\int_{0}^{r}\left(\frac{r^{m}}{s^{m-1}}-\frac{s^{m+1}}{r^{m}}\right)g\left(\beta_{2}g^{2}-\alpha+\beta^{'}f^{2}\right)ds,\quad r\in[0,\delta(C)).  \label{g1}
\end{equation}
And the local existence near 0 of \eqref{g1} is easy to get if we choose the initial function $g=Cr^{m}$, as $r\rightarrow 0$. It is interested in $C>0$. Define the following two disjoint sets as
\begin{align}
 S_{5}&\triangleq \left\{C>0|g' \ \text{becomes non-positive before }g\ \text{reaches}\ B\right\},\\
 S_{6}&\triangleq \left\{C>0|g\ \text{crosses}\ B\  \text{before}\ g'\ \text{becomes zero}\right\}.
\end{align}
Continuity in $C$ implies that $S_{5}$ and $S_{6}$ are open. Furthermore, they are
non-empty.
\par

\noindent $\bf{Step.2}$
We claim that $S_{5}$ contains small $C$.
Since $f\in C[0,\infty)$ and $\lim_{r\to\infty}f(r)=A$, set
\begin{equation}
\epsilon=-A+\sqrt{A^2+\frac{\beta_2 B^2}{4|\beta'|}}>0.
\end{equation}
Then there exists $R_{\epsilon}>1$ such that
\begin{equation}
A-\epsilon<f(r)<A+\epsilon,\quad r>R_{\epsilon}.
\end{equation}
Denote $R_1=k\max\{R_{\epsilon},\frac{2m}{\sqrt{\beta_1}B}\}$ where $k$ is an integer over $1$, and we will choose it later. For $C=0$, we know that $g\equiv 0$. The continuous dependence to the solution on the parameter allows us to choose sufficiently small $C$ such that
\begin{equation}
0\le |g|\le \frac{B}{4R_1},\quad for\ all\ r\ge 0.
\end{equation}
There exists $I^+$, such that
\begin{equation}
0\le g\le \frac{B}{4R_1},\quad r\in I^+.
\end{equation}
Recalling $\beta_2g^2(\infty)+\beta'f^2(\infty)=\beta_2B^2+\beta'A^2=\alpha$, we obtain
\begin{eqnarray}
g''+\frac{1}{r}g'&=&g\left(\frac{m^2}{r^2}+\beta_2 g^2-\alpha+\beta'f^2\right),\nonumber\\
&\le&g\left(\frac{m^2}{R_1^2}+ \frac{\beta_2B}{16R_1^2}^2-\alpha+\beta'A^2+|\beta'|\frac{\beta_2 B^2}{4|\beta'|}\right),\nonumber\\
&\le&g\left(\frac{\beta_2B^2}{4}+\frac{\beta_2B^2}{16R_1^2}-\beta_2B^2+\frac{\beta_2B^2}{4}\right),\nonumber\\
&\le&g\left(\frac{\beta_2B^2}{16R_1^2}-\frac{\beta_2B^2}{2}\right),\nonumber\\
&=&g\left(\frac{1-8R_1^2}{16R_1^2}\beta_2B^2\right)<0.\quad r>R_1.
\end{eqnarray}
Thus, we define the comparison function
\begin{equation}
y''+\frac{1}{r}y'+\frac{8R_1^2-1}{16R_1^2}\beta_2B^2y=0,
\end{equation}
which is a Bessel equation of order 0. As $r>\frac{R_1}{B\sqrt{(8R_1^2-1)\beta_2}}$, we have
\begin{equation}
y\sim \sqrt{\frac{8R_1}{\pi B\sqrt{(8R_1^2-1)\beta_2}r}}\cos\left(\frac{B}{4R_1}\sqrt{(8R_1^2-1)\beta_2}r-\frac{\pi}{4}\right)
\end{equation}
Define
\begin{align}
j_0=\inf \left\{j \Bigg|\frac{(8j+7)\pi R_1}{B\sqrt{(8R_1^2-1)\beta_2}}-R_1>0,\quad j=0,1,2\right\},
\end{align}
and let
\begin{equation}
z(j_0)=\frac{(8j_0+7)\pi R_1}{B\sqrt{(8R_1^2-1)\beta_2}}.
\end{equation}
Note that we can choose a proper $k$ such that
\begin{equation}
\sqrt{\frac{8R_1}{\pi Br\sqrt{(8R_1^2-1)\beta_2}}}> \frac{B}{4R_1}.
\end{equation}
We could find a suit interval $I_j=(z(j_0), z(j_{j_0+1}))\cap I^+$ as a comparison interval, in which we have $y(r)>g(r),\ r\in I_j$. However, $y(r)$ is the oscillating solution to the Bessel equation, and we conclude that there exists sufficiently small $C$, such that $g'(r)<0$ at some $r$ which implies $S_5$ is non-empty.

\medskip

\par

\noindent $\bf{Step.3}$
$S_{6}$ contains large $C$. Make the scaling $r=|C|^{-\frac{1}{m}}t$, so \eqref{g1} becomes
\begin{equation}
g(t)=t^{m}+\frac{1}{2m}\int_{0}^{t}
\left(\frac{t^{m}}{\tau^{m-1}}-\frac{\tau^{m+1}}{t^{m}}\right)\frac{1}{|C|^{\frac{2}{m}}}g\left(\beta_{2}g^{2}-\alpha+\beta^{'}f^{2}\right)d \tau.\label{g}
\end{equation}
We have $f \le Rr^{\lambda}\le R|C|^{-\frac{\lambda}{m}}t^{\lambda}$, for $r \le1$, which means $t\le |C|^{\frac{1}{m}}$
, we have
\begin{equation}
|C|^{-\frac{2}{m}}\beta'f^{2}\le |C|^{-\frac{2}{m}}\beta^{'}|C|^{-\frac{2\lambda}{m}}t^{2\lambda}R^{2}
\le \frac{\beta^{'}R^{2}}{|C|^{\frac{2}{m}}}\to 0,\
\text{when}\ C\to \infty.
\end{equation}
Similar to the proof of non-emptiness of $S_4$. If $|C|=\infty$, we can convert \eqref{g} to
\begin{equation}
g_{tt}+\frac{1}{t}g_{t}-\frac{m^2}{t^2}g =  |C|^{-\frac{2}{m}}g\left(\beta_{2}g^{2}-\alpha\right). \label{g2}
\end{equation}
We consider a formal power series expansion of $g$ near 0, $g=\sum_{i=k}^{\infty}a_it^i$, where $k$ needs to be determined and we get
\begin{align}
g&\sim t^m-\frac{|C|^{-\frac{2}{m}}\alpha}{4m+4}t^{m+2}+o(t^{m+2}),\nonumber\\
   &\sim |C|(r^{m}-\frac{\alpha}{4m+4}r^{m+2})+o(r^{m+2}).
\end{align}
Continuity guarantees that there exists some positive constant $N$, independent of the choice of $f$, such that ${\beta^{\prime}R^{2}}\mbox{/}{|C|^{\frac{2}{m}}}\le {1}\mbox{/}{N}$, then $|C|^{\frac{2}{m}}>N\beta^{\prime}R^{2}$. When $r$ sufficiently small, we conclude that there exists sufficiently large $C$, such that $g$ crosses $B$. One obtains there exists large $C$ in it and $S_{6}$ is non-empty.
\par

\noindent $\bf{Step.4}$
Since the connected set $C>0$ cannot consist of two open, disjoint and non-empty sets, there must exist some value of $C$ neither in $S_{5}$ nor $S_6$. We denote $C_{0}$ for this kind of value, which satisfies
\begin{equation}
0\le g<B,\quad
 g'>0.
\end{equation}
Because $g$ is increasing with an upper bound, then we use $g_{\infty}$ to denote the limit when $r\rightarrow\infty$. Let $r \to \infty$ at both sides of \eqref{f2} we have
\begin{equation}
\frac{1}{r}\left(rg_{r}\right)_{r}=g_{\infty}\left(\beta_2g_{\infty}^{2}-\alpha+\beta'A\right),
\end{equation}
where $g_{\infty}\in (0, B]$.
If $g_{\infty}\neq B$, then $\beta_2g_{\infty}^{2}-\alpha+\beta^{'}A\triangleq \tilde{a}<0$, and
\begin{equation}
\frac{1}{r}\left(rg_{r}\right)_{r}=\tilde{a},\quad \left(rg_{r}\right)_{r}\sim \tilde{a}r,\quad rg_{r}\sim \frac{\tilde{a}}{2}r^2,\quad g_{r}\sim \frac{\tilde{a}}{2}r,
\end{equation}
which is a contradiction for $g'>0$, so $g_{\infty}=B$.
\par

\noindent $\bf{Step.5}$
Finally, we want to prove the uniqueness of $g$. Suppose that there are two different solutions $g_1, g_2$ with $g_1(0)=g_2(0)=0$, $g_1(\infty)=g_2(\infty)=B$, and make the change $h=r^{\frac{1}{2}}g$ in \eqref{f2},
then we obtain
\begin{eqnarray}
r^2h_1''=h_1\left[\left(\beta'f^2-\alpha\right)r^2+\beta_2h_1^2r+m^2-\frac{1}{4}\right],  \label{h1}\\
r^2h_2''=h_2\left[\left(\beta'f^2-\alpha\right)r^2+\beta_2h_2^2r+m^2-\frac{1}{4}\right].  \label{h2}
\end{eqnarray}
Denote $H=h_1-h_2$, and then $H$ satisfies
\begin{equation}
r^2H''=H\left[\left(h_1^2+h_1h_2+h_2^2\right)\beta_2r+\left(\beta'f^2-\alpha\right)r^2+m^2-\frac{1}{4}\right]. \label{H}
\end{equation}
We know $H$ oscillates more slowly than $h_1$ from Sturm comparison theorem. However, $h_1$ has no zero when $r>0$, then $H$ cannot vanish at any finite point. Without loss of generality, one assumes that $H>0$ initially, which means $h_1>h_2\ (g_1>g_2)$. Multiplying $H$, $h_1$ at each sides of \eqref{h1} and \eqref{H} \cite{mcleod2001}, respectively, we get
\begin{equation}
(H'h_1-Hh_1')'=H''h_1-Hh_1''>0,
\end{equation}
which implies that $H'h_1-Hh_1'$ is increasing. Because $H(0)=h_1(0)=0$, then $h_1^{2}\left(\frac{H}{h_1}\right)'=H'h_1-Hh_1'\ge0$, which means ${H}\mbox{/}{h_1}$ is increasing which is equivalent to the fact that $\frac{g_1-g_2}{g_1}$ is increasing, so $\frac{g_1-g_2}{g_1}$ should have a strict positive limit, which is a contradiction as $r\rightarrow\infty$. So $g$ is unique, and then for any given $f$, $C_0$ is unique where $|C_0|\le (N\beta^{\prime})^{\frac{m}{2}}R^{m}$, proved Lemma 4.1.

\par

\noindent $\bf{Step.6}$
We give an estimate for $g$. Note that $\lim_{r\to 0}|r^{-m}g|=C_0$ and in any compact interval $[0,1]$, there exists $M>0$, such that $|r^{-m}g|\le M$, $r\in [0,1]$, where $M$ depends only on the parameters $(N.R, \alpha,\beta',\beta_1,\beta_2)$. Also, we can find a uniform upper bound for the term
\begin{equation}
|\beta_{2}g^{2}-\alpha+\beta^{'}f^{2}|\le \beta_{2}B^{2}+\alpha+|\beta'|R^2\triangleq \hat{M}(R),\quad r\in [0,1].
\end{equation}
Thus we obtain
{\setlength\arraycolsep{2pt}
\begin{eqnarray}
0<\left|g-C_0r^m\right|&\le &\frac{1}{2m}\int_{0}^r\left|\left(\frac{r^{m}}{s^{m-1}}-\frac{s^{m+1}}{r^{m}}\right)g\left(\beta_{2}g^{2}-\alpha+\beta^{'}f^{2}\right)\right|ds  \nonumber \\
   &=&\frac{1}{2m}\int_{0}^r\left|sr^m\left(1-\frac{s^{2m}}{r^{2m}}\right)\frac{g}{s^{m}}\left(\beta_{2}g^{2}-\alpha+\beta^{'}f^{2}\right)\right|ds \nonumber \\
  & \le & \frac{M\hat{M}(R)}{4m}r^{m+2}.\label{est1}
\end{eqnarray}}\qed

\medskip

\begin{lemma}\label{4.2}
For the $g(f)$ (denote $g$ for simplification below) in Lemma 4.1, we can find a unique $\tilde{f}\in C[0, \infty)$, satisfying
\begin{equation}
\tilde{f}''+\frac{1}{r}\tilde{f}'-\frac{n^2}{r^2}\tilde{f}=\tilde{f}\left[\beta_1(\tilde{f}^2-1)+\beta'g^2\right],   \label{tilf11}
\end{equation}
and
\begin{equation}
\tilde{f}(0)=0, \quad \tilde{f}(\infty)=A, \quad 0\leq \tilde{f}<A, \quad \tilde{f}'>0,
\end{equation}
with $r^{-n}\tilde{f}$ decreasing. If we choose a proper large $R$ in Lemma 4.1, then we get $\tilde{f}\le R^{*}(R)r^{n}$, $r\le 1$.
\end{lemma}

\medskip

\noindent$\bf{Proof.}$ We divide it into six steps.
\par

\noindent $\bf{Step.1}$For the homogeneous form of \eqref{tilf11}, its solutions are $r^n, r^{-n}$, we can convert it formally into the integral
equation
\begin{equation}
\tilde{f}=Dr^{n}+\frac{1}{2n}\int_{0}^{r}
\left(\frac{r^{n}}{s^{n-1}}-\frac{s^{n+1}}{r^{n}}\right)\tilde{f}\left[\beta_1\left(\tilde{f}^2-1\right)+\beta'g^{2}\right]ds.\quad 0\le r\le \delta(D). \label{tilf2}
\end{equation}
Equation \eqref{tilf2} can be solved by iteration, at least for $r$ sufficiently small where $\tilde{f}$ has a natural estimate
\begin{equation}
\tilde{f}=Dr^{n}+O\left(r^{n+2}\right).\label{tilf0}
\end{equation}
We are only interested in the case $D>0$. Define the sets
{\setlength\arraycolsep{2pt}
\begin{align*}
& S_{7}\triangleq \{D>0|\tilde{f}'\ \text{reaches} \ 0\ \text{before}\ \tilde{f}\ \text{reaches}\ A\},\\
& S_{8}\triangleq \{D>0|\tilde{f}\ \text{crosses}\ A\ \text{before}\ \tilde{f}'\ \text{reaches}\ 0\}.
\end{align*}}
The solution depends continuously on the parameter $D$, clearly they are both open sets. From the structure of $S_7$, $S_8$, we know $S_{7}\cap S_{8}=\varnothing$. Then we claim that $S_7$ and $S_8$ are both non-empty.

\par

\noindent $\bf{Step.2}$ $S_7$ is  non-empty. Let $\tilde{f}=Du$ in \eqref{tilf1}, then $u$ satisfies
\begin{equation}
u''+\frac{1}{r}u'-\frac{n^2}{r^2}u = u\left[\beta_1(D^2u^2-1)+\beta'g^2\right]. \label{u1}
\end{equation}
Let $D\to 0$, then the \eqref{u1} becomes
\begin{equation*}
u''+\frac{1}{r}u'-\frac{n^2}{r^2}u = u\left(\beta'g^2-\beta_1\right).
\end{equation*}
We interest in $u>0$. If $u\le 0$, $S_3$ is non-empty obviously. Then we get
{\setlength\arraycolsep{2pt}
\begin{eqnarray}
\begin{split}
u''+\frac{1}{r}u'+\left(\beta_1-\beta'B^2-\frac{n^2}{r^2}\right)&u=u\beta_1(g^2-B^2)<0, \label{u2} \\
u\sim r^n, \quad r&\to 0.
\end{split}
\end{eqnarray}}
Consider
{\setlength\arraycolsep{2pt}
\begin{align}
\begin{split}
v''+\frac{1}{r}v'+&\left(\beta_1-\beta'B^2-\frac{n^2}{r^2}\right)v=0, \label{u3} \\
&v\sim r^n, \quad r\to 0.
\end{split}
\end{align}}
Assume that $v>0$ in $[0,R_D)$. Comparing \eqref{u2} with \eqref{u3}, we can infer that $v>u$ in $(0, \delta)$ when $\delta$ is sufficiently small. In fact if there exists $r_1$, such that $v(r_1)=u(r_1)$, then one constructs
\begin{equation}
c=\sup\{\hat{c}>0|\ \hat{c}v\le u,\ r\in (0, r_1)\},
\end{equation}
which implies that there exists $r_2$, such that
\begin{equation}
cv=u,\quad cv'=u',\quad cv''\le u'',\quad r=r_2.
\end{equation}
 But it is impossible,
{\setlength\arraycolsep{2pt}
\begin{eqnarray}
0&=&\left.c\left[v''+\frac{1}{r}v'+\left(\beta_1-\beta'B^2-\frac{n^2}{r^2}\right)v\right]\right|_{r=r_2},\nonumber\\
&\le&\left.u''+\frac{1}{r}u'+\left(\beta_1-\beta'B^2-\frac{n^2}{r^2}\right)u\right|_{r=r_2},\nonumber\\
&=&u(r_2)\beta'\left(g^2(r_2)-B^2\right)<0.
\end{eqnarray}}
A contradiction! So $v>u>0$, $r\in[0,R_2)$, where $R_2<R_D$.  Moreover, noting that $\beta_1-\beta'B^2=\frac{\beta_1(\beta_1\beta_2-\beta'\alpha)}{\beta_1\beta_2-\beta'^2}>0$, then $v$ is a solution of standard Bessel equation \cite{chen1994}, so there must exist $r_3$,
such that $u'(r_3)=0$, which means $S_7$ is non-empty from openness when $D$ is small.\par

\par

\noindent $\bf{Step.3}$
Next we prove $S_8$ is non-empty. Let $r=D^{-\frac{1}{n}}t$ in \eqref{tilf2}, then
\begin{equation}
\tilde{f}(t)=t^{n}+\frac{1}{2n}\int_{0}^{t}\left(\frac{t^{n}}{\tau^{n-1}}-\frac{\tau^{n+1}}{t^{n}}\right)D^{-\frac{2}{n}}\tilde{f}\left[\beta_1(\tilde{f}^2-1)+\beta'g^{2}\right]d \tau. \label{tilft}
\end{equation}
Similarly, if $|D|=\infty$, we can convert \eqref{tilft} to
\begin{equation}
\tilde{f}_{tt}+\frac{1}{t}\tilde{f}_{t}-\frac{n^2}{t^2}\tilde{f}=D^{-\frac{2}{n}}\beta_{1}\tilde{f}(\tilde{f}^{2}-1).
\end{equation}
When $t$ is small, make a series expansion of $\tilde{f}$. Then the indeterminate coefficients we get are
\begin{align}
\tilde{f}&\sim t^n-\frac{D^{-\frac{2}{n}}\beta_1}{4n+4}t^{n+2}+o\left(t^{n+2}\right),\nonumber\\
         &\sim |D|\left(r^{n}-\frac{\beta_{1}}{4n+4} r^{n+2}\right)+o\left(r^{n+2}\right).
\end{align}
So we can find a sufficiently large $D$, such that $\tilde{f}$ crosses $A$ when $r$ sufficiently small and $S_8$ is non-empty.
\par

\noindent $\bf{Step.4}$
If the solutions satisfy the boundary conditions, then $D_0\notin S_7\cup S_8$, and $|D_{0}|^{\frac{2}{n}}<N\beta^{\prime}B^{2}$ which means
\begin{equation}
 0\le \tilde{f}<A,\quad \tilde{f}'>0.
\end{equation}
$\tilde{f} $ is increasing with an upper bound. We define $\lim_{r \rightarrow \infty}\tilde{f}(r)=\tilde{f}_{\infty}$, then $\tilde{f}_{\infty}\in (0, A]$. Let $r\to \infty$ in the two sides of equation \eqref{tilf1}, we obtain
\begin{equation}
\frac{1}{r}\left(r\tilde{f}_{r}\right)_{r}=\tilde{f}_{\infty}\left[\beta_1(\tilde{f}_{\infty}^{2}-1)+\beta'B^2\right].
\end{equation}
After same argument we obtain $\tilde{f}_{\infty}=A$.
\par

\noindent $\bf{Step.5}$
It is required to prove the uniqueness of $\tilde{f}$. Assume there exists $\tilde{f}_1,\ \tilde{f}_2$  satisfying
$\tilde{f}_1(0)=\tilde{f}_2(0)=0,\  \tilde{f}_1(\infty)=\tilde{f}_2(\infty)=A$. Let $q=r^{\frac{1}{2}}\tilde{f}$ in \eqref{tilf1}, then it converts to
\begin{equation}
r^2q''=q\left[r\beta_1(q^2-r)+\beta'g^2r^2+n^2-\frac{1}{4}\right].\label{q}
\end{equation}
Denote $Q=q_1-q_2$, then
\begin{equation}
r^2Q''=Q\left[r\beta_1\left(q_1^2+q_1q_2+q_2^2-r\right)+\beta'g^2r^2+n^2-\frac{1}{4}\right];\label{Q}
\end{equation}
$q_1$ satisfies
\begin{equation}
r^2q_1''=q_1\left[r\beta_1(q_1^2-r)+\beta'g^2r^2+n^2-\frac{1}{4}\right].\label{q1}
\end{equation}
We see: $Q$ oscillates more slowly than $q_1$, and $q_1$ has only one zero at $r=0$, so $Q$ can
not equal to 0 at finite point from Sturm comparison theorem. Without loss of generality, we assume $Q>0$ initially. Comparing \eqref{Q} and \eqref{q1}, we know that ${Q}\mbox{/}{q_1}$ is increasing monotonically, then $\frac{\tilde{f}_1-\tilde{f}_2}{\tilde{f}_2}$ should has a strictly positive limit,
which is a contradiction and one gives the uniqueness of $\tilde{f}$.
\par

\noindent $\bf{Step.6}$
We do some preparation before giving the estimate for $\tilde{f}$, as $r\le 1$. From \eqref{tilfn}, we obtain
\begin{align}
\frac{\tilde{f}}{r^{n}}&=D_0+\frac{1}{2n}\int_{0}^{r} \left(\frac{1}{s^{n-1}}-\frac{s^{n+1}}{r^{2n}}\right)\tilde{f}\left[\beta_1(\tilde{f}^{2}-1)+\beta'g^{2}\right]ds, \nonumber\\
&=D_0+\frac{1}{2n}\int_{0}^{r}\left(s-\frac{s^{2n+1}}{r^{2n}}\right)\frac{\tilde{f}}{s^{n}}\left[\beta_1(\tilde{f}^{2}-1)+\beta'g^{2}\right]ds. \label{tilfr}
\end{align}
Derive \eqref{tilfr}, one gets
\begin{align}
\left(\frac{\tilde{f}}{r^{n}}\right)'=\int_{0}^{r} \frac{s^{n+1}}{r^{2n+1}}\tilde{f}\left[\beta_1(\tilde{f}^{2}-1)+\beta'g^{2}\right]ds,
\end{align}
and from above we conclude that $\tilde{f}\mbox{/}{r^{n}}$ is decreasing because of sign of the forcing term. However, $\lim_{r\to 0} r^{-n}\tilde{f}=D_0$, thus we know $r^{-n}\tilde{f}\le D_0$, $r\le 1$. Also, term $|\beta_1(\tilde{f}^2-1)+\beta'g^2|$ has an upper bound which is denoted by $\bar{M}$ and it is a function of $\hat{M}(R)$, i.e. for simplification: $\bar{M}(R)$. Recalling \eqref{est1}, \eqref{tilf2} and the fact we just got, one has
{\setlength\arraycolsep{2pt}
\begin{eqnarray}
\left|\tilde{f}-D_{0}r^{n}\right|&=&\frac{1}{2n}\left|\int_{0}^{r}\left(\frac{r^{n}}{s^{n-1}}-\frac{s^{n+1}}{r^{n}}\right)\tilde{f}\left[\beta_1(\tilde{f}^2-1)+\beta'g^2\right]ds\right|\nonumber\\
&\le& \frac{1}{2n}\int_{0}^{r}r^{n}s\frac{\tilde{f}}{s^{n}}\left|\beta_1(\tilde{f}^2-1)+\beta'g^2\right|ds  \nonumber\\
 &\le& \frac{D_0\beta_1\bar{M}}{2n}r^{n}\int_{0}^{r} s ds \nonumber\\
 &=& \frac{D_0\beta_1\bar{M}}{4n}r^{n+2},\quad 0\le r\le 1.
\end{eqnarray}
It means
\begin{equation}
|r^{-n}\tilde{f}|\le D_0(1+\frac{\beta_1\bar{M}(R)}{4n}r^2),\quad 0\le r\le 1.
\end{equation}
and Lemma 4.2 is proved.\qed

\medskip

\begin{remark}\label{4.2}
There is a map from $R$ to $R^*$, that is
\begin{align}
\sigma:\quad &\mathbb{R}\rightarrow \mathbb{R},\nonumber\\
&R\mapsto R^*.
\end{align}
We conclude that the map $\sigma$ always has at least one fix point recorded by $R^*$ if we choose a suitable $r_0\in [0,1]$, and then we get $\sigma(R^*)=R^*$.
\end{remark}
\medskip

\noindent $\bf{Proof \ of\  Theorem\ 2.2.}$ Define the $Banach$ space $\mathscr{B}$
\begin{equation}
\mathscr{B}=\big\{f\in C[0, \infty)\ \big| \ \exists M>0,\ r^{-\lambda}(1+r^{\lambda})f\le M\big\},\label{banach}
\end{equation}
with the norm
\begin{equation}
\|f\|_{\mathscr{B}}=\sup_{r\in [0, \infty)}\big\{\big|r^{-\lambda}(1+r^{\lambda})f\big|\big\}.\label{norm}
\end{equation}
One chooses a non-empty bounded closed convex set $S$ as
\begin{equation}
S=\big\{f\in B:\big|r^{-n}f\big|\le R^{*}, \quad r\le r_0; \quad f(\infty)=A,\quad f< A, \quad r^{-n}f\ \text{is decreasing}\big\}.\label{set}
\end{equation}
Here $R^{*}$ in $S$ is chosen uniformly from the \emph{Remark 4.1}. Evidently, $S$ is non-empty, bounded, closed and convex. We need to prove the following\par
($i$) The mapping $\Phi$ maps $S$ into itself, and $\Phi:f\mapsto \tilde{f}$.\par
($ii$) $\Phi$ is continue.\par
($iii$) $\Phi$ is precompact.\par
Then the Schauder fixed point theorem ensures that $\Phi$ has at least one fixed point, and the theorem is proved.

Part($i$) is guaranteed by Lemma 4.1 and Lemma 4.2. Similarly, the continuity is concluded from the ODE theory. To prove part($iii$), given any sequence $\{f_k\}\in S$, which is bounded in the norm of $\mathscr{B}$, then we see that $\tilde{f}_k'$ is bounded in any finite interval $(0,r)$, $(r\neq 0)$. From the Differential mean value theorem, one knows that $\tilde{f}_{k}$ is equicontinuous and by the Ascoli-Arzela theorem, there exists a subsequence of $\{\tilde{f}_k\}$, converging uniformly to some $\tilde{f}$ in any compact subinterval of $(0, \infty)$, and the difficulties occur \emph{at} $r\to 0$ and $r\to \infty$.

When $r\to 0$, given $\epsilon >0$, there exists $\hat{r}$ sufficiently small, such that
\begin{equation}
\sup_{r\le r_0}|r^{-\lambda}(\tilde{f_k}-\tilde{f})|=\sup_{r\le r_0}|r^{-n}(\tilde{f_k}-\tilde{f})r^{n-\lambda}|\le \hat{r}^{n-\lambda}
R^{*}<\epsilon.
\end{equation}
When $r\to \infty$, rewrite \eqref{tilf11} to consider the asymptotic property of $\tilde{f}$
\begin{equation}
\tilde{f}_k''+\frac{1}{r}\tilde{f}_k'-\frac{n^2}{r^2}\tilde{f}_k=\tilde{f}_k\left[\beta_1(\tilde{f}_k^2-A^{2})+\beta'(g^2-B^{2})\right].\label{tilf3}
\end{equation}
To acquire the fact we would like to get the estimate for $\tilde{f}_k$ as $r\to \infty$, let
\begin{equation}
e(r)=\beta_1 \tilde{f}_k(\tilde{f}_k+A),
\end{equation}
and consider comparison function
\begin{equation}
F=\mu r^{-\sigma},\quad r>0,\quad \mu>0.\label{F0}
\end{equation}
We get
\begin{equation}
F''+\frac{1}{r}F'-\frac{n^{2}}{r^{2}}F=\frac{\sigma^{2}-n^{2}}{r^{2}}F.\label{F1}
\end{equation}
Combining \eqref{tilf3}, one gets
\begin{eqnarray}
(F+\tilde{f}_k-A)''+\frac{1}{r}(F+\tilde{f}_k-A)'&=&\left[\frac{\sigma^{2}-n^{2}}{r^{2}}+\frac{n^{2}}{\mu}r^{\sigma -2}-e(r)\right
]F+\beta' \tilde{f}_k(g^{2}-B^{2}) \nonumber \\
&&+(e(r)+\frac{n^{2}}{r^{2}})(F+\tilde{f}-A).
\end{eqnarray}
Set $\sigma=2$, because of the fact that $0\le g<B$, there exists $\bar{r}>1$, for any $r> \bar{r}$, we can choose $\mu$ such that
\begin{equation}
\frac{4-n^{2}}{r^{2}}+\frac{n^{2}}{\mu}-e(r)<0.\ \
\end{equation}
Choose $\mu$ large in \eqref{F0} such that $(F+f-A)(\bar{r})>0$ and we infer that
\begin{equation}
(F+\tilde{f}_k-A)''+\frac{1}{r}(F+\tilde{f}_k-A)'\le (e(r)+\frac{n^{2}}{r^{2}})(F+\tilde{f}_k-A).
\end{equation}
From maximum principle, we obtain
\begin{equation}
(F+\tilde{f}_k-A)(r)\ge 0,\quad r\ge \bar{r},
\end{equation}
 then
\begin{equation}
-\mu r^{-2}\le \tilde{f}_k(r)-A \le 0.
\end{equation}
Also, we derive that $g\ge B-\nu r^{-2}$, $r>\bar{r}$ in the similar approach. For the upper bound, we can choose $\hat{\mu}$ small, such that \begin{equation}
(F+\tilde{f}_k-A)(\bar{r})<0,
\end{equation}
 and
\begin{equation}
\left[\frac{4-n^{2}}{r^{2}}+\frac{n^{2}}{\hat{\mu}}-e(r)\right]F+\beta' \tilde{f}_k(g^{2}-B^{2})\ge 0,\quad r>\bar{r}.
\end{equation}
Using Maximum principle we get
\begin{equation}
\tilde{f}\le A-\hat{\mu} r^{-2}.
\end{equation}
Thus it is right when $r\to\infty$ and one gets
\begin{equation}
\tilde{f}_k=A+O(r^{-2}).\label{ef}
\end{equation}
We could get the asymptotic property of $g$ the same as $f$. By using the asymptotic property, we can compute the quantization by the way we did in the previous section and prove Theorem 2.2.

\end{document}